# Determining Mills' Constant and a Note on Honaker's Problem

Chris K. Caldwell

Department of Mathematics and Statistics

University of Tennessee at Martin

Martin, TN 38238

USA

caldwell@utm.edu

Yuanyou Furui Cheng

Department of Mathematics,

Temple University,

Philadelphia, PA 19122.

cfy1721@gmail.com

**Abstract**

In 1947 Mills proved that there exists a constant $A$ such that $\lfloor A^{3^n} \rfloor$ is a prime for every positive integer $n$. Determining $A$ requires determining an effective Hoheisel type result on the primes in short intervals—though most books ignore this difficulty. Under the Riemann Hypothesis, we show that there exists at least one prime between every pair of consecutive cubes and determine (given RH) that the least possible value of Mills' constant $A$ does begin with $1.3063778838$. We calculate this value to $6850$ decimal places by determining the associated primes to over $6000$ digits and probable primes (PRPs) to over $60000$ digits. We also apply the Cramér-Granville Conjecture to Honaker's problem in a related context.

## 1 Introduction

In 1947 Mills [19] proved that there exists a constant $A$ such that $\lfloor A^{3^n} \rfloor$ is a prime for every positive integer $n$. Mills' proof used Ingham's result [17] that there is always a prime in the interval $(x, x + kx^{5/8})$ for some constant $k$. Since $k$ was not determined explicitly, Mills could not determine $A$ or even a range for $A$. Many authors citing Mills' constant either follow his example and remain silent about $A$'s value (e.g., [11], [24], [31]), or explicitly state that it is unknown (e.g., [1], [5]). But recently a few authors have begun to state that $A$ is approximately 1.3063778838 (e.g., [8], [10], [28], [33] and Sloane A051254). They do this despite the fact that it is currently impossible to compute $A$ without any further unproven assumptions, and despite the fact that for each $c \geq 2.106$, there are uncountably many possible values $A$ for which $\lfloor A^{c^n} \rfloor$ is prime for each positive integer $n$.

The authors who approximate Mills' constant $A$ also implicitly define *Mills' constant* to be the least $A$ such that $\lfloor A^{3^n} \rfloor$ is prime for all positive integers $n$. Mills' original article contained no numerics—it only showed that such an $A$ exists.

In Section 2 we observe that, assuming the Riemann Hypothesis, it easily follows that there is at least one prime between $x^3$ and $(x+1)^3$ for every $x \geq 0.26$ (Lemma 5). We then use this to determine Mills' constant $A$ to over 6850 decimal places in Section 3. This requires a substantial computing effort to calculate the associated *Mills' primes* $\lfloor A^{3^n} \rfloor$ explicitly for $n = 1, 2 \ldots, 10$. In particular we find prove the following:

**Theorem 1.** *Assume the Riemann Hypothesis. The minimum Mills' constant (for the exponent c=3) begins with the following* 600 *digits:*

$$
\begin{array}{lllll}
1.3063778838 & 6308069046 & 8614492602 & 6057129167 & 8458515671 \\
3644368053 & 7599664340 & 5376682659 & 8821501403 & 7011973957 \\
0729696093 & 8103086882 & 2388614478 & 1635348688 & 7133922146 \\
1943534578 & 7110033188 & 1405093575 & 3558319326 & 4801721383 \\
2361522359 & 0622186016 & 1085667905 & 7215197976 & 0951619929 \\
5279707992 & 5631721527 & 8412371307 & 6584911245 & 6317518426 \\
3310565215 & 3513186684 & 1550790793 & 7238592335 & 2208421842 \\
0405320517 & 6890260257 & 9344300869 & 5290636205 & 6989687262 \\
1227499787 & 6664385157 & 6619143877 & 2844982077 & 5905648255 \\
6091500412 & 3788524793 & 6260880466 & 8815406437 & 4425340131 \\
0736114409 & 4137650364 & 3793012676 & 7211713103 & 0265228386 \\
6154666880 & 4874760951 & 4410790754 & 0698417260 & 3473107746
\end{array}
$$

We also find the two probable-primes which are (most likely) the next terms in the sequence of Mills' primes. The last of these has 61684 digits.

Concerning other problems that involve prime gaps, let $p_0 = 2$, $p_n$ be the $n$-th odd prime number and $g_n = p_{n+1} - p_n$, the gap between consecutive prime numbers. The Cramér-Granville conjecture [16] states that $g_n \leq M \log^2 n$ for some constant $M > 1$. The last section is an application of this conjecture on prime gaps to Honaker's problem.

Honaker's problem [9] is to find all trios of consecutive prime numbers $p < q < r$, such that $p | (qr + 1)$. We call such triples of consecutive prime numbers *Honaker trios* and conjecture that there are only three Honaker trios: $(2, 3, 5)$, $(3, 5, 7)$, and $(61, 67, 71)$. We establish that there are only three Honaker trios $(p, q, r)$ for $p \leq 2 \times 10^{17}$ and prove:

**Theorem 2.** *The Cramér-Granville conjecture with any constant $M$ implies that there are only a finite number of Honaker trios. If also $M \leq 199262$, then there are exactly three trios.*

## 2 Gaps Between Prime Numbers

Let $s(q)$ represent the next prime after $q$: $s(q) = \min\{p \text{ prime} \mid p > q\}$. Define the sequence of maximal gap primes as follows:

$$q_j = \begin{cases} 2, & \text{if } j = 0; \\ 3, & \text{if } j = 1; \\ \min\{p \text{ prime} \mid s(p) - p > s(q_{j-1}) - q_{j-1}\}, & \text{if } j > 1. \end{cases}$$

In other words $\{q_k\}$ is the sequence of primes that are followed by maximal gaps (Sloan A002386). This sequences starts 2, 3, 7, 23, 89, 113, 523, 887, 1129, 1327 ... and exists for all subscripts $k$ because $\limsup_{p \to +\infty} s(p) - p = +\infty$ [15]. The list of maximal gap primes $q_j$ has been extended through all primes below $2 \times 10^{17}$ by Herzog and Silva [30] with much of their work verified by Nicely and others [22, 23].

With these definitions it is clear that

$$\frac{s(p) - p}{\log^2 p} \leq \frac{s(q_k) - q_k}{\log^2 q_k},$$

for every $q_k \leq p < q_{k+1}$. Now using Nicely's tables of maximal gaps [22] we can easily verify the following result:

**Lemma 3.** *For $11 \leq p_n < 2 \times 10^{17}$, we have $\frac{s(p)-p}{\log^2 p} \leq 0.92064$.*

To explicitly determine Mills' constant we must know that there is always a prime between each pair of successive cubes from some established point onward. This follows easily from a result of Schoenfeld [29]. Recall

$$\text{li}(x) := \lim_{\epsilon \to 0^+} \left( \int_0^{1-\epsilon} \frac{dt}{\log t} + \int_{1+\epsilon}^x \frac{dt}{\log t} \right).$$

**Lemma 4** (Schoenfeld). *Assume the Riemann Hypothesis. For $x \geq 2657$, we have*

$$\text{li}(x) - \frac{\sqrt{x}\log x}{8\pi} < \pi(x) < \text{li}(x) + \frac{\sqrt{x}\log x}{8\pi}.$$

(Much stronger bounds are now available [27], but this is sufficient for our current needs.)

Using this Lemma we see that if $x > 2657^{1/3}$, then

$$\pi((x+1)^3) - \pi(x^3) > \text{li}((x+1)^3) - \text{li}(x^3) - \frac{3}{4\pi}(x+1)^{3/2}\log(x+1)$$
$$\geq \int_{x^3}^{(x+1)^3} \frac{dt}{\log t} - \frac{3}{4\pi}(x+1)^{3/2}\log(x+1)$$
$$\geq \frac{3x^2 + 3x + 1}{3\log(x)} - \frac{3}{4\pi}(x+1)^{3/2}\log(x+1).$$

This last lower bound is an increasing function of $x$ and it is greater than one for $x > 2657^{1/3}$. After using a computer program to check the smaller values of $x$ we have shown the following:

**Lemma 5.** *Assume the Riemann Hypothesis. There is at least one prime between $x^3$ and $(x+1)^3$ for every $x \geq x_0 = 2^{1/3} - 1$.*

This lemma states that, under the Riemann Hypothesis, there are prime numbers between consecutive cubes; we will see that this is necessary to calculate Mills' constant in the next section.

Without the Riemann Hypothesis it is still possible to use explicit upper bounds on the zeta function [6, 13] to get a version of Lemma 5 with a far larger value of $x_0$, roughly $10^{6000000000000000000}$ [7, 26]. In the next section we will show that calculating Mills' constant involves explicitly finding a prime larger than the cube-root of this bound. This bound so dramatically exceeds current computing abilities, that any current calculation of Mills' constant must involve an unproven assumption. This bound $x_0$ should get much smaller as the bounds on the zeta function are improved.

## 3 Mills' constant

We begin this section with a lemma.

**Lemma 6.** *If $x > 1$ and $c > 2$, then $1 + x^c + x^{c-1} < (1+x)^c$.*

*Proof.* Dividing by $x^c$ and replacing $x$ with $1/x$ we arrive at the equivalent inequality $0 < (1+x)^c - (1+x+x^c)$ $(0 < x < 1)$. The inequality clearly holds when $c = 2$ (because it reduces to $x > 0$) and when $x = 0$. Now if $x > 0$, differentiate the right side with respect to $c$ to get $(1+x)^c \log(1+x) - x^c \log(x)$, which is clearly positive, so the inequality above holds for all $c > 2$. □

It will be useful to next recall a proof of a simple generalization of Mills' theorem.

**Theorem 7.** *Let $S = \{a_n\}$ be any sequence of integers satisfying the following property: there exist real numbers $x_0$ and $w$ with $0 < w < 1$, for which the open interval $(x, x+x^w)$ contains an element of $S$ for all real numbers $x > x_0$. Then for every real number $c > \min\left(\frac{1}{1-w}, 2\right)$ there is a number $A$ for which $\lfloor A^{c^n} \rfloor$ is a subsequence of $S$.*

*Proof.* (We follow Ellison & Ellison [12].) Define a subsequence $b_n$ of $S$ recursively by

(a) $b_1$ is equal to the least member of $S$ for which $b_1^c$ is greater than $x_0$.

(b) $b_{n+1}$ is the least member of $S$ satisfying $b_n^c < b_{n+1} < b_n^c + b_n^{wc}$.

Because $c > 1/(1-w)$ and $c > 2$, (b) can be written as:

$$b_n^c < b_{n+1} < 1 + b_{n+1} < 1 + b_n^c + b_n^{wc} < 1 + b_n^c + b_n^{c-1} < (1+b_n)^c$$

the last inequality following from Lemma 6. For all positive integers $n$ we can raise this to the $c^{-(n+1)}$th power to get

$$b_n^{c^{-n}} < b_{n+1}^{c^{-(n+1)}} < (1+b_{n+1})^{c^{-(n+1)}} < (1+b_n)^{c^{-n}}.$$

This shows that the sequence $\{b_n^{c^{-n}}\}$ is monotonic and bounded, therefore converges. Call its limit $A$. Finally,
$$b_n < A^{c^n} < b_n + 1,$$
so $\lfloor A^{c^n} \rfloor = b_n$, yielding the chosen subsequence of $S$ and completing the proof. □

When constructing the sequence $\{b_n\}$ in the previous proof, the condition (a) can be relaxed in two important ways. First, by replacing the words 'the least' with 'any', we see there are infinitely many choices for $b_1$, and hence for the resulting value $A$. Second, it is not necessary that $b_1^c > x_0$ be satisfied by $b_1$, as long as the terms satisfying (b) exist to a term $b_n$ which does satisfy $b_n^c > x_0$. This will someday be important in removing the assumption of the Riemann Hypothesis from Theorem 1 by calculating a sequence of primes $\{b_i\}_{i=1}^n$ that extends to the cube root of the bound $x_0$ discussed at the end of the previous section.

It is proved by Baker *et. al.* [3] that $p_{n+1} - p_n = O(p_n^{0.525})$. From the above theorem, one obtains the following proposition.

**Proposition 8.** *For every $c \geq 2.106$, there exist infinitely many $A$'s such that $\lfloor A^{c^n} \rfloor$ is a prime for every $n$.*

Wright [32] showed that the set of possible values of the constants $A$ and $c$ in this proposition (and several generalizations [25], [28], [31]) have the same cardinality as the continuum and are nowhere dense. So authors approximating "Mills' constant" must first decide how to choose just one. Mills specified the exponent $c = 3$. Mills also used the lower bound $k^8$ for the first prime in the sequence, where $k$ is the *integer* constant in Ingham's result [17]. So to follow his proof literally would require that the first prime be 257 (and producing a constant $A \approx 6.357861928837$). But all authors offering an approximation agree implicitly on starting with the prime 2 and then choosing the least possible prime at each step, so we take this as the definition of *Mills' constant*.

The sequence of minimal primes satisfying the criteria of the proof with $c = 3$ (Sloan A051254) begins with

$$\begin{aligned}
b_1 &= 2, \\
b_2 &= 11, \\
b_3 &= 1361, \\
b_4 &= 2521008887, \\
b_5 &= 16022236204009818131831320183, \\
b_6 &= 4113101149215104800030529537915953170486137\backslash \\
&\quad 96235397599331359499488277040407483
2568499.
\end{aligned}$$

These first six terms were well known. We have now shown that assuming the Riemann Hypothesis, prime numbers exist between consecutive cubes, so (with RH) we know this sequence can be continued indefinitely.

To make these fast growing primes $b_n$ easier to present, define a sequence $a_n$ by $b_{n+1} = b_n^3 + a_n$. The sequence $a_n$ begins with $3, 30, 6, 80, 12, 450, 894, 3636, 70756$ (Sloan A108739). The primality of the new terms $b_7$, $b_8$ and $b_9$ (2285 digits) were proven using the program

Titanix [18] which is based on an elliptic curve test of Atkin ([2], [20]). The test for $b_9$, which at the time was the third largest proven 'general' prime, was completed by Bouk de Water in 2000 using approximately five weeks of CPU time. The certificate was then verified using Jim Fougeron's program Cert-Val.

François Morain has verified the primality of the next term $b_{10}$ (6854 digits, July 2005) using his current implementation of fastECPP [14, 21]. The computation was done on a cluster of six Xeon biprocessors at 2.6 GHz. Cumulated CPU time was approximately 68 days (56 for the DOWNRUN, 12 for the proving part).

In late 2004 Phil Carmody used his self-optimizing sieve generator to generate an appropriate sieve. Using it and pfgw he verified the above results (as PRPs) and found the next two probable-primes in the sequence $b_n$. They are the values corresponding to $a_{10} = 97220$ and $a_{11} = 66768$. These yield probable primes of 20562 and 61684 digits respectively, so their primality may not be proven for some time.

The first ten terms of $\{b_n\}$ are sufficient to determine Mills' constant $A$ to over 6850 decimal places because we have from the proof of Theorem 7:

$$b_n^{c^{-n}} < A < (b_n + 1)^{c^{-n}}.$$

A quick calculation now completes the proof of Theorem 1.

## 4 Honaker's Problem

Let $g(p)$ be the length of the prime gap after $p$: $g(p) = s(p) - p$.

Cramér's original conjecture says that $g(p) = (1 + o(1)) \log^2 p$. This conjecture could be too strong to be valid, so Granville [16] proposed there is a bound $M$ for which $g(p) < M \log^2 p$ (we refer to this conjecture as the Cramér-Granville conjecture). As we discussed in Section 2, $M = 0.92064$ is valid for $11 \le p < 2 \times 10^{17}$. However, the constant involved in this conjecture may be $M = 2e^{-\gamma} \approx 1.123$ [16].

*Proof of Theorem 2.* Searching for Honaker trios with $p < 50$ we find only $(2, 3, 5)$ and $(3, 5, 7)$. In what follows we will assume $(p, q, r)$ is a Honaker trio and $p > 50$.

Reformulate the problem by letting $q = p + 2k$ and $r = q + 2l$. Since $p|(qr+1)$, we have $p|(4k^2 + 4kl + 1)$. Thus,

$$p \le 4k^2 + 4kl + 1. \tag{1}$$

If $k \ge l$, then $p \le 8k^2 + 1$ and $g(p) = 2k \ge \sqrt{(p-1)/2}$. Otherwise $k < l$, $p \le 4(l-1)^2 + 4(l-1)l + 1$ and again $g(q) = 2l \ge \sqrt{(p-1)/2}$. These gaps are bounded by the Cramér-Granville conjecture, so either $p$ or $q$ must satisfy:

$$M \ge \sqrt{\frac{p-1}{2}} \frac{1}{\log^2 p}. \tag{2}$$

The function on the right is increasing for $p > 50$, and is unbounded, so for any fixed $M$ the number of Honaker trios is finite.

Next suppose $p < 2 \times 10^{17}$. By Lemma 3 we know $M = 0.92064$ will suffice in this range, so equation (1) gives $p < 14142$. A computer search in this range adds the third trio $(61, 67, 71)$. Thus there are exactly three trios with $p < 2 \times 10^{17}$.

Finally, since any additional solution must satisfy $p > 2 \times 10^{17}$, equation 2 shows $M > 199262$. This completes the proof of Theorem 3. □

Clearly this simple proof can be extended to make similar statements about finding $k$ consecutive primes for which one of the $k$ divides the product of the other $k-1$ plus or minus a fixed integer.

Return to